\documentclass{article}

\usepackage{url} 
\usepackage{lineno}
\usepackage{soul}

\usepackage{authblk}

\usepackage{lmodern,lscape,textcomp,placeins}
\usepackage[T1]{fontenc}
\usepackage{multirow}

\usepackage{amsmath,amssymb,dsfont,units,booktabs}
\usepackage{tabularx}
\usepackage{appendix}
\usepackage{graphicx, color}
\usepackage{color}
\usepackage[utf8]{inputenc} 
\usepackage[round]{natbib}
\usepackage{algorithm}
\usepackage{algpseudocode}


\newcommand{\vt}{\mathbf{v}}

\newcommand{\taut}{\pmb{\tau}}

\newcommand{\sigmat}{\pmb{\sigma}}

\newcommand{\epsilont}{\pmb{\dot \epsilon}}

\newcommand{\te}{\pmb{e}}

\renewcommand{\div}{\operatorname{div}}
\newcommand{\tr}{\operatorname{tr}}

\title{Simulating sea-ice deformation in viscous-plastic sea-ice models with CD-grids}

\date{}

\author[1]{C. Mehlmann  \thanks{Corresponding author: carolin.mehlmann@ovgu.de}}
\author[2]{ G.Capodaglio}
\author[3,4]{S. Danilov}

\affil[1]{\footnotesize Otto-von-Guericke Universit\"at, Magdeburg, Germany}
\affil[2]{\footnotesize
Computational Physics and Methods Group, Los Alamos National Laboratory, Los Alamos, U.S.A.}
\affil[3]{\footnotesize Alfred Wegener Institute, Helmholtz Centre for Polar and Marine Research, Bremerhaven, Germany}
\affil[4]{Jacobs University, Bremen, Germany}

\begin{document}

\maketitle

\begin{abstract}%
Linear Kinematic Features (LKFs) are found everywhere in the Arctic sea-ice cover. They are strongly localized deformations often associated with the formation of leads and pressure ridges. Viscous-plastic sea-ice models start to produce LKFs at high spatial grid resolution, typically with a grid spacing below 5 km. A recent study showed that the placement of the variables on the grid plays an important role for the number of simulated LKFs. The study found that a nonconforming finite element discretization with a CD-grid placement (CD1) resolves more LKFs per degree of 
freedom compared to more common A,B and C-grids.
A new CD-grid formulation (CD2) has just been proposed based on a conforming subgrid discretization. 
To analyze the resolution properties of the new CD2 approach we evaluate runs from different models (e.g FESOM, MPAS) on a benchmark problem using  quadrilateral,  hexagonal and triangular meshes. We found that the CD1 setup simulates more deformation structure than the CD2 approximation. This highlights the importance of the type of spatial discretization for the simulation of LKFs.
Due to the higher number of degrees of freedom both CD-grids resolve more LKFs than traditional A,B and C-grids at fixed mesh level. This is an appealing feature as  high spatial mesh resolution is needed in viscous-plastic sea-ice models to simulate LKFs.
\end{abstract}

\section{Introduction}
Sea ice plays a crucial role in the climate system, as it acts as a buffer between the ocean and the atmosphere, influencing the exchange of heat, moisture and momentum \citep{Stroeve2018}.
Linear kinematic features (LKFs), such as leads and ridges, are build during drift in response to wind and ocean forcing.
These LKFs are important indicators of sea-ice deformation which are closely linked to the mechanical properties of sea ice. However, simulating these features in viscous-plastic sea-ice models has proven to be a significant challenge. 

In most climate models sea ice is characterized as a viscous-plastic two dimensional continuum, which is  represented either in the classical viscous-plastic (VP) \citep{Hibler1979} formulation or by the elastic-viscous-plastic (EVP) model modification \citep{Hunke1997}.  The viscous-plastic sea-ice model has been criticized in the last years 
for using assumptions that have no observational evidence \citep{Coon2007,Rampal2008,Feltham2008}. Alternative rheologies have also been proposed \citep{Dansereau2016, Rampal2016,Tsamados2013}. Even so, most practical applications are still using the (E)VP formulation and will continue to apply it in the foreseeable future \citep{Blockley2020}.

It has been demonstrated in different studies that (E)VP models are able to simulate aspects of observed LKFs \citep{Hutter2020}: once the resolution is high enough these models are able to reproduce observed deformations and spatial-temporal scaling laws \citep{Bouchat2017,Bouchat2022}.
However the simulated LKFs depend on numerical details used in the model realization  such as the solver convergence \citep{Koldunov,Lemieux2012}, the pressure paramterization \citep{Hutchings2005}, the mesh resolution \citep{Wang2014} or the placement of the sea-ice variables on the grid \citep{Mehlmannetal2021}. 

 This paper is a follow up study of \cite{Mehlmannetal2021}, where the authors observed that CD-grids, which place the velocity variables at the edge midpoint, were the most effective in simulating LKFs across varying grid spacing and degrees of freedom. We found that the statement made in \cite{Mehlmannetal2021} is not complete: in addition to the edge placement, the spatial discretization also plays a crucial role for the amount of simulated deformation. Our observations hold true for different types of grids, including quadrilateral, triangular, and hexagonal meshes. 

In the context of the viscous-plastic sea-ice model, an increase of the mesh resolution leads to a larger amount of simulated deformation. However, instead of increasing the resolution, more elaborated spatial discretizations can be an efficient alternative. \cite{Mehlmannetal2021} found that CD-grids resolve more deformation than approaches which represent the velocity by its normal component at the edge midpoint (C-grids) or setups with velocities located on vertices (A-grids, B-grids). The study also analyzed the simulated deformation on triangular meshes and found that the triangular CD-grid resolves more LKFs than the triangular A-grid. The effect can mainly be attributed to the higher number of degrees of freedom (dof).

Additionally, the authors observed that the CD-grid discretization has appealing resolution properties \citep{Mehlmannetal2021}: the  quadrilateral CD-grid simulates a similar amount of deformation structure than quadrilateral B-grid and C-grid approaches, but on a grid twice as coarse. The same is true for the triangular CD-grid compared to the triangular A-grid.

The study of \cite{Mehlmannetal2021} considered a CD-grid based on a nonconforming finite element discretization (Crouzeix-Raviart element) \citep{Mehlmann2021}. We refer to this approximation as CD1 in the following.  A new CD-grid formulation has been proposed by \cite{Capodaglio2022}, we will refer to this approach as CD2. It differs from CD1 in the spatial discretization of strain rates and stresses. The approach applies a subgrid discretization based on Wachpress functions \citep{Dasgupta2003} or piecewise linear bases \citep{BAILEY2008}. The question arises on how the resolution properties of CD2 approaches compare to those of CD1 and traditional A-grids, B-grids and C-grids.

We extend the study of \cite{Mehlmannetal2021} and evaluate the ability of the new CD2 approximations to simulate LKFs.
The numerical analysis is performed on a benchmark problem using quadrilateral, triangular and hexagonal meshes. The CD2 simulations on quadrilateral and hexagonal grids
 are obtained with the sea-ice module of the \emph{Model for Prediction Across Scales} \citep[MPAS-Seaice]{Turner2021,Capodaglio2022}. In case of triangular meshes, we compute the CD2 approximation with the \emph{Finite-Volume Sea Ice–Ocean model} \citep[FESOM]{Danilov2015}. Please note that the CD2 discretization in FESOM \citep{Danilov2023} is inspired by the development of \cite{Capodaglio2022} and differs in numerical details in the subgrid discretization. 
 
 The CD2 simulations are compared to runs conducted with the \emph{Los Alamos Sea Ice Model} \citep[CICE]{Hunke2015}, the sea-ice module of the \emph{Icosahedral Nonhydrostatic Weather and Climate Model} \citep[ICON]{Mehlmann2021},  the sea-ice module of the \emph{Massachusetts Institute of Technology
general circulation model} \citep[MITgcm]{LOSCH2010}, and the setup realized in the academic software library \emph{Gascoigne} \citep{Gascoigne}. Note that
MPAS-Seaice relies on the EVP method. FESOM and ICON solve the VP formulation with a modified EVP solver while all the other models are based on a implicit VP formulation. Our analysis shows that, for the simulation of LKFs, the differences between VP and EVP are minor compared to the chosen discretization of the velocity. 

This paper is structured as follows. Section \ref{sec:model} introduces the viscous-plastic sea-ice model and outlines the used numerical discretizations and methods applied for the analysis. 
Section \ref{sec:num} presents a numerical analysis of the data. A discussion is given in Section \ref{sec:diss}. The paper ends with a conclusion in Section \ref{sec:conc}.

\section{Governing equations}\label{sec:model}
We consider a simplified sea-ice model, where sea-ice is characterized by three variables: sea-ice velocity $\vt$, sea-ice thickness $H$ and sea-ice concentration $A$. The sea-ice dynamics is described by the following system of equations
\begin{linenomath}
  \begin{align}
    m \partial_t \vt 
    + f_c \te_z\times \vt&=
    \div\,\sigmat + F,
   \label{eq:mom}\\
    \partial_t A + \div\,(\vt A) &= 0, \quad  A\le 1,\label{eq:A}  \\
    \partial_t H + \div\,(\vt H) & = 0. \label{eq:h}
  \end{align}
\end{linenomath}
By setting the right-hand side to zero in  equations (\ref{eq:A}) and (\ref{eq:h}), all thermodynamic source terms are neglected.
The ice mass per unit area is $m=\rho_\text{ice}H$, where $\rho_{\text{ice}}=\unit[900]{kg/m^{3}}$ is the density. For the other terms in the momentum equation, $f_c$ is the Coriolis
parameter, $\te_z$ is the vertical ($z$-direction) unit vector, and $\sigmat$ is the internal stress. The external forces are collected in
\[F= A\,\taut(\vt)  -\rho_\text{ice}H g\nabla \tilde H_g,\]
where $g$ is the gravitational acceleration, $\tilde H_g$ is the sea surface height and $\taut(\vt)$ describes the oceanic and atmospheric stresses. The internal stresses $\sigmat$ are related to the strain rates \begin{equation}\label{eq:strain}
\epsilont=\frac{1}{2}(\nabla \vt +\nabla \vt^T)
\end{equation} by the  viscous-plastic (VP) material law \citep{Hibler1979}
\begin{linenomath}\label{eq:stress}
  \begin{align}
 \sigmat=2 \eta\epsilont +(\zeta-\eta)\tr(\epsilont)I-\frac{P}{2}I, \quad P=\frac{P_0\Delta}{2(\Delta+\Delta_{min})}.
  \end{align}
\end{linenomath}
Note that the superscript $T$ in (\ref{eq:strain}) indicates the transpose,  $I$ is the identity matrix and
$P$ is the replacement pressure, which has been introduced by \cite{HiblerIp95}  to avoid divergence of ice in the absence of forces.
The viscosities $\zeta,\eta$ are given by 
\begin{linenomath}
\begin{equation}\label{eq:visc}
  \zeta=\frac{P_0}{2 (\Delta^2+\Delta^2_{min})^{\frac{1}{2}}},\qquad  \eta=e^{-2}\zeta,\qquad
  P_0(H,A)= P^\star H \exp\big(-C(1-A)\big),
\end{equation}
\end{linenomath}
where $e=2$ is the ratio of the elliptic yield curve, $P^\star=27500$ N/m$^2$ is the ice strength parameter and $C=20$. The parameter $\Delta_{min}=2\times 10^{-9}$ is the viscous limit of the plastic regime, and 
\begin{linenomath}
\begin{equation}
\Delta^2= (\epsilont_{11}^2+\epsilont_{22}^2)(1+e^{-2}) + 
4\epsilont_{12}^2\, e^{-2}+2\epsilont_{11}\epsilont_{22}(1-e^{-2}). 
\end{equation}
\end{linenomath}

The elastic-viscous-plastic formulation \citep{Hunke1997} has been introduced to regularize the VP rheology. 
\begin{equation}
    \partial_t \sigmat +\frac{e^2}{2T_{evp}}\sigmat+\frac{1-e^2}{4T_{evp}}\tr(\sigmat)I+\frac{P}{4T_{evp}}I=\frac{\zeta}{T_{evp}} \epsilont,
\end{equation}
where $\tr(\cdot)$ denotes the trace and $T_{evp}$ is the relaxation time that determines the transition time from the elastic regime to the VP rheology.
The viscous-plastic material law is recovered for $\partial_t\sigmat=0$.

\subsection{Method section}
We consider the benchmark problem introduced by \cite{Mehlmannetal2021}: the idealized test case models the initial stage of sea-ice deformation. In this setup, a domain of size 512 km $\times$ 512 km is covered with a thin layer of sea ice. A cyclone moves diagonally through the domain, sea-ice deforms and multiple LKFs are built during the 2 days of simulations. We compare the shear deformation \begin{linenomath}\begin{align}\label{shear}
     \epsilont_\text{shear}=\sqrt{(\epsilont_{11}-\epsilont_{22})^2 +4 \epsilont_{12}^2},
\end{align}\end{linenomath} simulated by different model setups after these two days. The shear deformation 
is analyzed either visually or by using a detection algorithm outlined by \cite{Hutter2019}.

The configuration of the algorithm used here interpolates the model data on a 2 km regular grid and detects LKFs by using image recognition tools in Python. 
After detection, the number and total length of LKFs are provided. The detection algorithm only identifies LKFs that are wider than one pixel to avoid the detection of noise. More details on configuration of the detection algorithm can be found in \citep{Mehlmannetal2021}.

We will compare the simulation in terms of two aspects: the amount of deformation with respect to the mesh size and the resolved structure in terms of degrees of freedom (dof). The simulated deformation with respect to the mesh is given by the total number or total length of LKFs for a given grid resolution. This metric highlights which discretization resolves more deformation structure on a fixed mesh. The second metric, the  resolved structure per dof, is given by the number or length of LKFs for a given amount of dof. 
This measure qualitatively characterizes the numerical effort used for the simulation of LKFs.


\subsection{Discretization}
The most common approach to solve the coupled sea-ice system (\ref{eq:mom})-(\ref{eq:h}) is to split the equations in time. First the approximation of the sea-ice momentum equation (\ref{eq:mom}) is computed, followed by solving the transport equations (\ref{eq:A})-(\ref{eq:h}). Due to stability concerns, fully explicit time stepping methods for the momentum equation are avoided as an extremely small time step is necessary in this case \citep{Hibler2000}. There are currently two ways to address this issue. One of them relies on using an implicit time discretization and iterative methods such as Picard solvers \citep{Lemieux2009,Hibler1991} and Newton methods \citep{Lemieux2010,Losch2014,MehlmannRichter2016newton,Shih2022}.

The other one still uses an explicit discretization, but relies on the EVP model \citep{Hunke1997}, in which an artificial elastic term is added to the VP rheology, allowing for an explicit discretization of the momentum equation with relatively large time steps. However, the original EVP model does not simulate the same deformation as VP models \citep{Lemieux2012,Boullion2013}. Therefore, a modified version of the EVP method, referred to as the mEVP solver, was developed to ensure convergence to the solution of the VP model \citep{Boullion2013,Kimmritz2015,Lemieux2012}.

The mEVP method is designed in such a way that numerical stability and convergence are addressed separately. In the original EVP formulation of \cite{Hunke1997}  numerical stability needs to be ensured by taking a sufficiently high number of subcycles. By slightly modifying the parameters choice of the EVP setup stable solutions can be produced with reduced iteration count \citep{Danilov2021}. \cite{Danilov2021} showed that (m)EVP methods lead to qualitative similar approximations of the benchmark problem if a sufficiently high number of subcycles is used. We extend the work and demonstrate that this is also true if the benchmark problem with a viscous-plastic rheology is solved implicitly. Using coarse meshes (27 km resolution) \cite{KIMMRITZ2017} showed that VP and mEVP approximations lead to quantitatively similar results in realistic Arctic setup.

We use the MPAS-Seaice in the B-grid \citep{Turner2021} and CD-grid \citep{Capodaglio2022} configuration on quadrilateral and hexagonal meshes. In both setups, PWL basis functions \citep{BAILEY2008} are adopted for the computation of the divergence of the stress in the momentum equation (\ref{eq:mom}). The sea-ice drift is calculated based on the EVP method using 500 subcycles. As suggested by \cite{Danilov2022} we select $T_{evp}=25$ min in the EVP algorithm. Please note that by increasing $T_{evp}$ it will be possible to further decrease the number of subcycles.
For the advection, an incremental remapping scheme is used \citep{Turner2021}. Due to the current lack of an incremental remapping scheme for a CD-grid, in such a case the velocities are first interpolated from the edges to the vertices and then the standard B-grid advection scheme is applied.

We compare the MPAS runs to simulations conducted with implicit solvers (CICE, Gascoigne, MITgcm) or the mEVP method (FESOM, ICON).  In the case of the mEVP approximation, we apply 100 subcycles per time step. As in the setup of \cite{Mehlmannetal2021} we solve the benchmark problem with a 2 minute time step. 
Further details on the different model configurations can be found in \citep{Mehlmannetal2021}. A description of the CD2 setup in FESOM is provided by \cite{Danilov2023}. 


\subsection{Structured and unstructured grids}\label{sec:mesh}
\begin{figure}
\begin{center}
\includegraphics[scale=0.4]{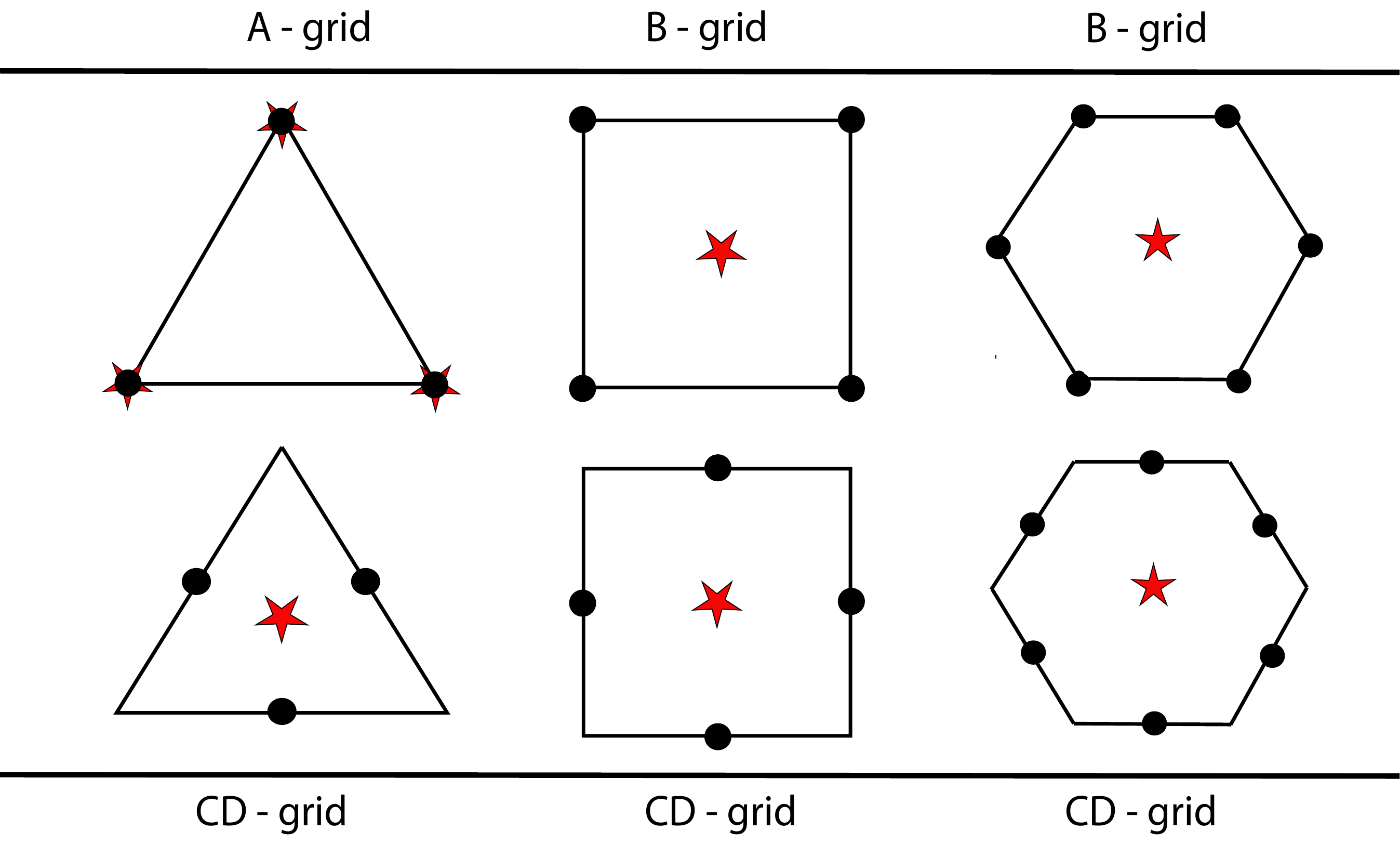}
\caption{Different staggering explored in this manuscript. We indicate the placement of the latitude/zonal  velocity $\vt=(u,v)$ and the staggering of the tracers by $\bullet$, and {\color{red}$\star$}, respectively.\label{fig:grids}}   
\end{center}
\end{figure}
We evaluate the two CD-grid discretizations on both structured quadrilateral grids and unstructured meshes. An overview of the setups considered by us is given in Figure \ref{fig:grids}. 
Traditionally sea-ice models have been formulated on structured quadrilateral meshes with either an Arakawa B-grid or C-grid discretization, which are used for example in CICE \citep{Hunke2015} or the MITgcm \citep{LOSCH2010}, respectively. Recent developments have also taken in account unstructured meshes made out of triangles or hexagons.  For example, unstructured triangular grids are applied in the sea-ice component of FESOM \citep{Danilov2015,Danilov2023} with an A, B and CD-grid approximation or in the sea-ice module of ICON \citep{Mehlmann2021,MehlmannGutjahr2022} with a CD-grid discretization. Hexagonal meshes are currently applied in MPAS-Seaice \citep{Turner2021,Capodaglio2022} based on a B-grid or a CD-grid discretization. 

The choice of grids in this study follows the setup used in \cite{Mehlmannetal2021}. In the quadrilateral case we analyze the benchmark problem on grids with 8 km, 4 km, and 2 km side, which correspond to 4096, 16384, and 65536 cells, respectively.
The CD-grid has $\frac{8N}{2}$ dof and the B-grid contains of $\frac{8N}{4}$ dof, where $N$ is the number of cells. The dof are calculated as follows. We have $8$  velocity components per cell; this amount is divided by two in the CD-grid setup as two neighboring cells share a dof. In the B-grid case a dof connects four adjacent cells. 

On triangular grids we use  meshes with a side length of 8 km, 4 km and 2 km with 9490, 37926, 151630 cells, respectively. The triangular CD-grid consists of $\frac{6N}{2}$ dof, the B-grid has $2N$ dof.

The runs on regular hexagonal grids are performed on meshes  with a distance between hexagon centers of 8 km, 4 km and 2 km. These meshes are made out of 4464, 18396, 74676 cells, respectively. The hexagonal B-grids and CD-grids have $12N/3$ dof and $12N/2$ dof.

We will compare CD-grid discretizations on triangular and  hexagonal meshes which have approximately the same amount of dof.
In the framework of MPAS-Seaice we also consider a B-grid mesh with the same dof as the CD-grid. This configuration is achieved by doubling  the number of cells in the case of quadrilateral meshes and increasing it by a factor of $1.5$ in case of hexagonal elements.

\section{Numerical evaluation}\label{sec:num}
\begin{figure}[tp]
  \begin{center}
  \includegraphics[scale=1.0]{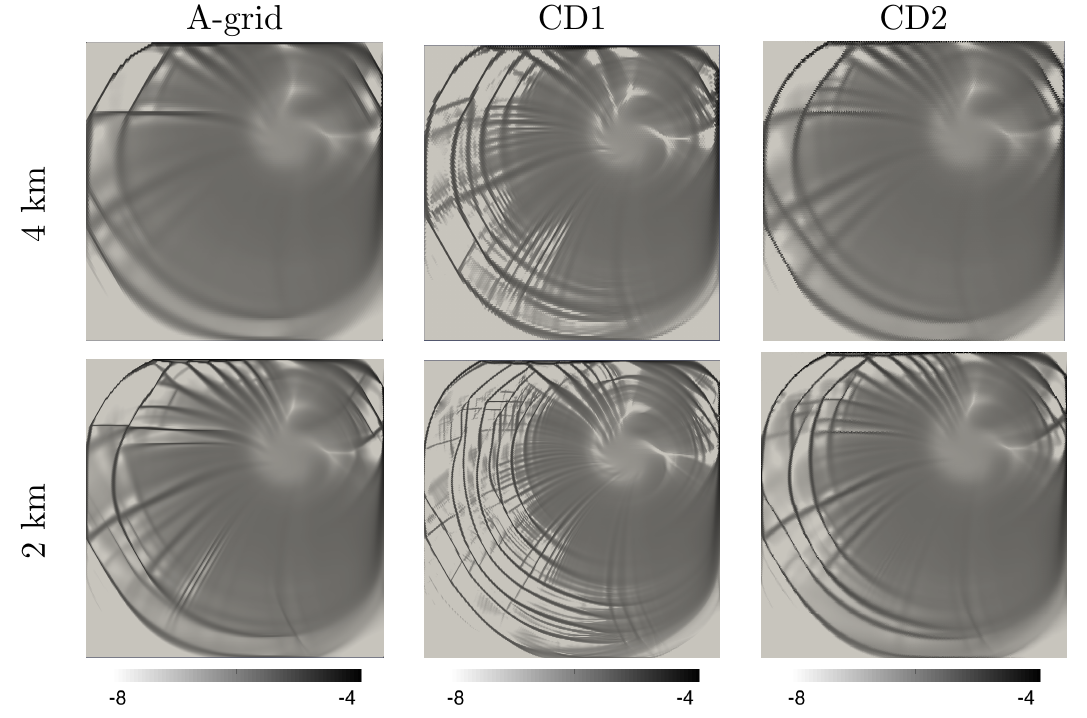}
 \caption{ {\small Shear deformation on a triangular mesh in FESOM. The shear deformation is given in $s^{-1}$ and plotted in logarithmic scale.
 \label{fig:trishear}}}
 \end{center}
\end{figure}
In this section we evaluate the ability of the CD1 and CD2 discretizations to reproduce LKFs. The analysis is conducted on structured and unstructured meshes. 
\subsection{ Triangular meshes}
We start by comparing the CD1 and CD2 discretizations in FESOM. Both setups place the velocity at the edge midpoint, use an upwind scheme to discretize (\ref{eq:A})-(\ref{eq:h}), and the iterative mEVP method to solve the momentum equation (\ref{eq:mom}). 
\begin{figure}
    \centering
    \includegraphics[scale=1.0]{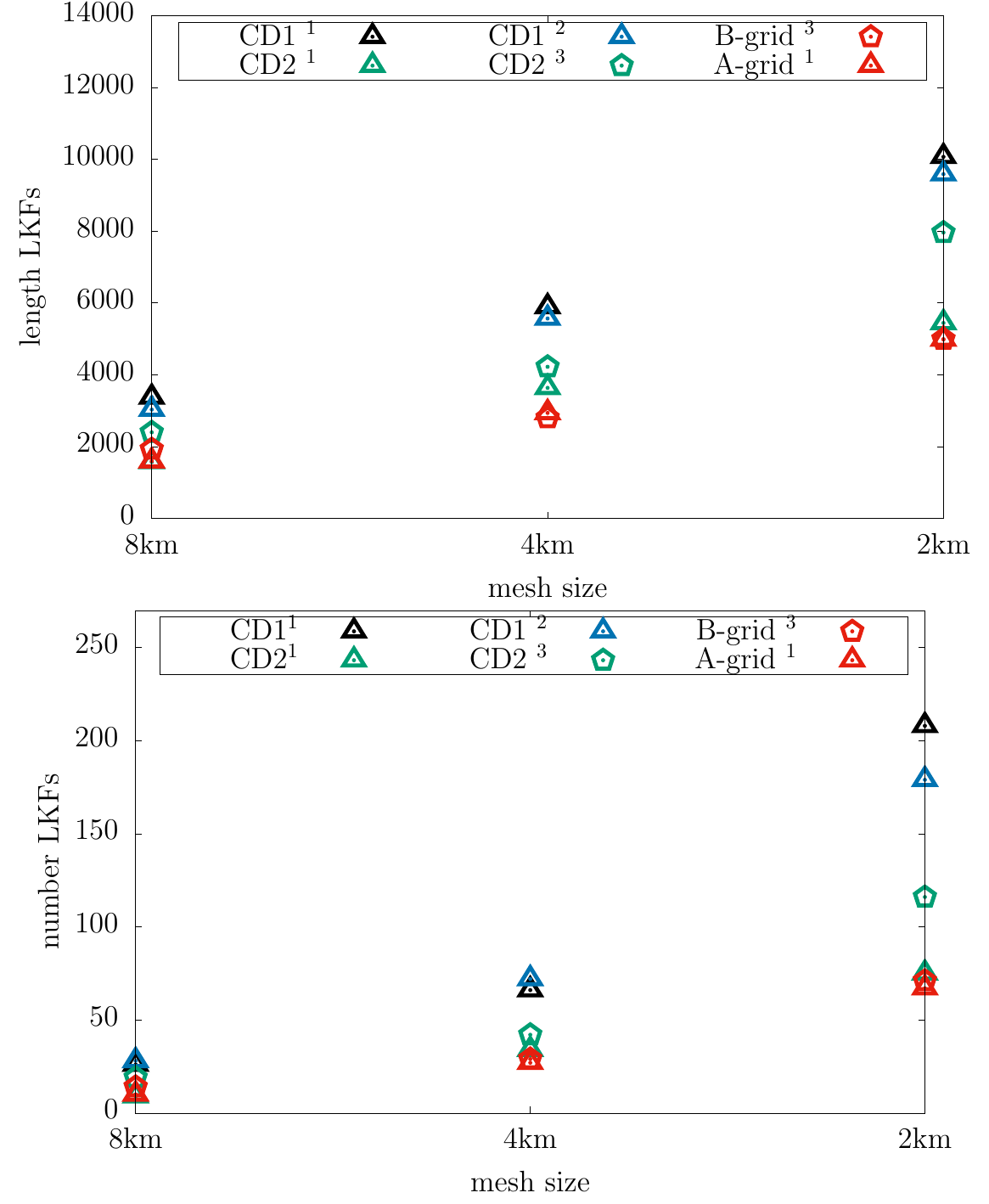}
    \caption{Total length and number of LKFs detected on triangular and hexagonal grids. Note that the CD-grid approaches have the same dof on hexagonal and triangular grids. The superscripts 1, 2, and 3 refer to the simulations carried out in the framework of FESOM, ICON and MPAS, respectively.}
   \label{fig:tri_LKFs}
\end{figure}
Figure \ref{fig:trishear} presents the simulated shear deformation of the two CD-grid approaches. The corresponding number of detected LKFs and their total length are given in Figure \ref{fig:tri_LKFs}.
As the CD1 and CD2 approximation  only differs in the handling of the spatial approximation, we
attribute the significant changes in the simulated deformation to the used spatial discretization.

Generally the largest number of LKFs is simulated with the CD1 setup followed by the CD2 approximation and the A-grid discretization. This observation is supported by the LKF detection algorithm. Overall the CD1 approach simulates finer structure than the CD2 and A-grid setup. The influence of the tracer placement to the amount of simulated LKFs is analyzed in \citep{Danilov2023,MehlmannDanilov2022}.
While in the CD2 setup the changes of the tracer location have a small influence on the amount of simulated LKFs, moving the tracer from cell to vertex smoothed out the LKFs in the CD1 configuration \citep{Danilov2023,MehlmannDanilov2022}.

\subsection{ Hexagonal meshes}
\begin{figure}[tp]
  \begin{center}
    \includegraphics[scale=1.0]{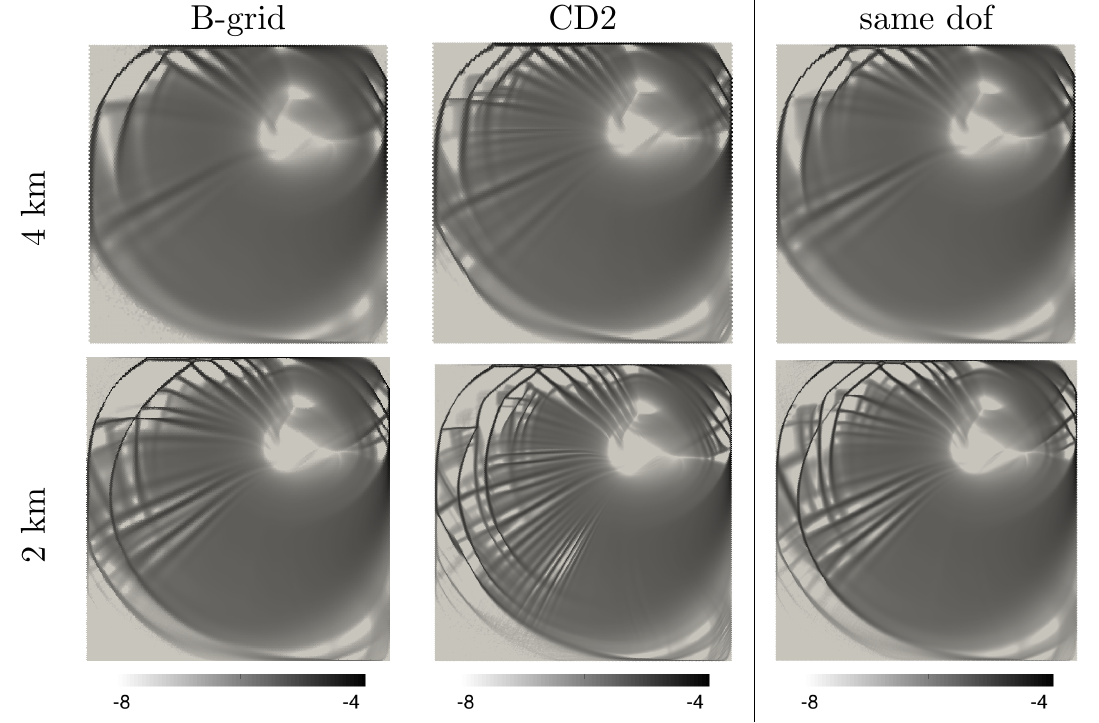}
 \caption{ {\small Shear deformation calculated on a hexagonal mesh in the framework of MPAS. The heading 'same dof' refers to a B-grid discretization with the same velocity dof as the CD-grid.  The shear deformation is given in $s^{-1}$ and plotted in logarithmic scale.
 \label{fig:hexshear}}}
 \end{center}
\end{figure}
\begin{figure}
    \centering
  \includegraphics[scale=1.0]{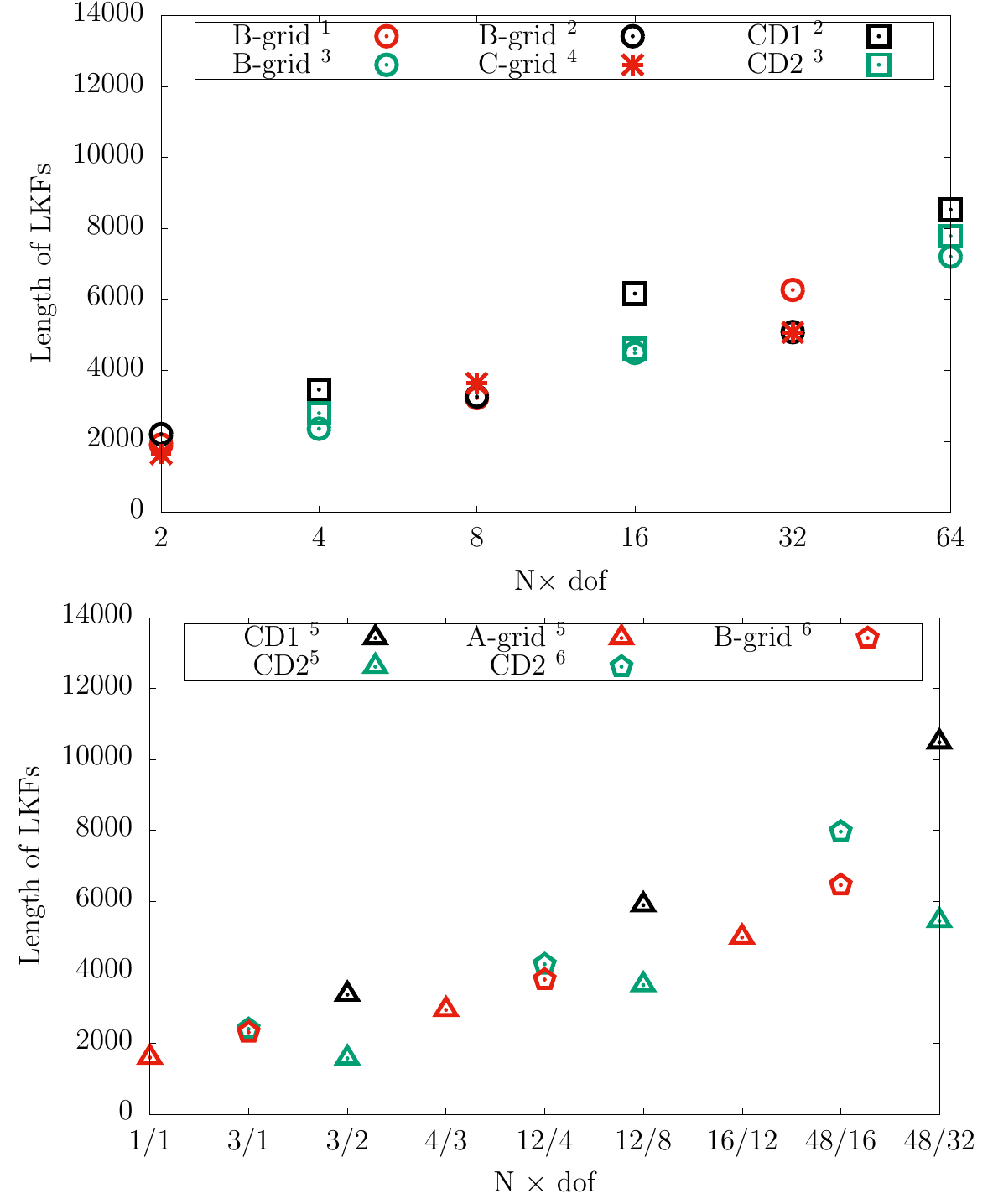}
   \caption{Detected number and total length of LKFs with respect to the degrees of freedom (dof).  N is the number of cells of the 8 km mesh. In the quadrilateral case (upper row) the numbers refer to the velocity dof. Note that all considered quadrilateral setups have the same tracer dof. In the triangular/hexagonal setup (bottom row) the numbers indicate the dof in the velocity/tracer components. The superscripts 1, 2, 3 and 4 in the upper plot refer to the simulations carried out in the framework of CICE, Gascoigne, MPAS and MITgcm respectively. The superscript 5 and 6 in the lower plot refer to the runs performed in FESOM and MPAS respectively.}
   \label{fig:quads_dof_LKFs}
\end{figure}
The B-grid and CD-grid approaches in MPAS only differ by the placement of the velocities and corresponding spatial discretization.  Figure \ref{fig:hexshear} shows that the MPAS CD2 discretization produces more structure than the B-grid setup. This is also reflected in the detected number and total length of LKFs presented in  Figure \ref{fig:tri_LKFs}.
By comparing the B-grid and the CD-grid approach on meshes with the same number of dof we see that the CD-grid setup still produces more structure than the B-grid approximation. This finding is supported by the result of detection algorithm presented in  Figure \ref{fig:quads_dof_LKFs}.

The MPAS CD2 implementation of the benchmark problem differs from the FESOM CD2 approximation by the choice of grid decomposition, the numerical realization of the subgrid discretization,
the used pseudo time-stepping (mEVP vs EVP) and the advection scheme. To compare result of CD-grids on hexagonal and triangular meshes we use grids with the same dof (see Section \ref{sec:mesh} for more details on the mesh choice).
The CD2 version in MPAS simulates less structure than the the CD1 setup on triangles. However it reproduces more deformation structure than the CD2  case on triangles (see Figure \ref{fig:hexshear}). This observation is also reflected in the output of the detection algorithm ( see Figure \ref{fig:tri_LKFs}). This difference is likely due to the different implementations of the CD2 discretization in MPAS and FESOM.

\subsection{Quadrilateral meshes}

\begin{figure}
  \begin{center}
 \includegraphics[scale=1.0]{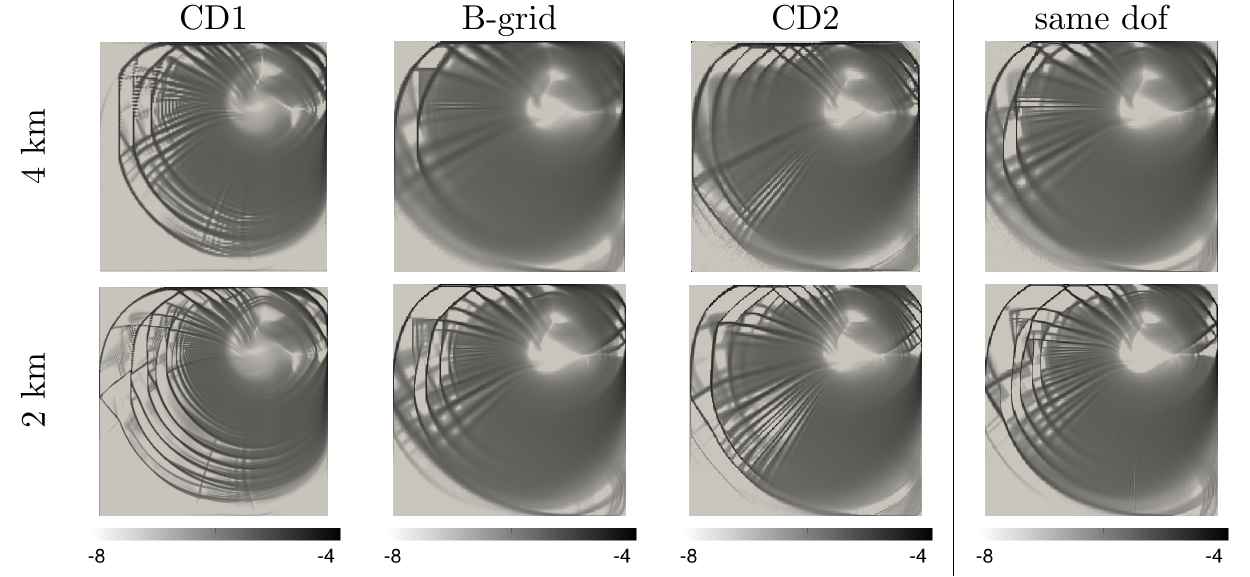}
  \end{center}
  \begin{center}
   \caption{ Shear deformation on quadrilateral meshes. The CD1 simulation is performed in Gascoigne whereas the CD2 and B-grid runs are carried out in the framework of MPAS. The heading "same dof" refers to a B-grid discretization with the same velocity dof as the CD-grids.  The shear deformation is given in $s^{-1}$ and plotted in logarithmic scale.  
   \label{fig:quads}}
  \end{center}
\end{figure}
We consider the B-grid and CD2 framework in MPAS and the CD1 setup in Gascoigne. The Gascoigne CD1 configuration is based on a non-conformal finite element discretization \citep{Mehlmannetal2021}, while the CD2 approximation in MPAS uses a conformal sub-grid discretization for the approximation of the stress and strain rates. The MPAS B-grid discretization differs from the MPAS CD2 setup only by the placement of the velocity degrees of freedom and the associated spatial discretization.

\begin{figure}
  \begin{center}
    \includegraphics[scale=1.0]{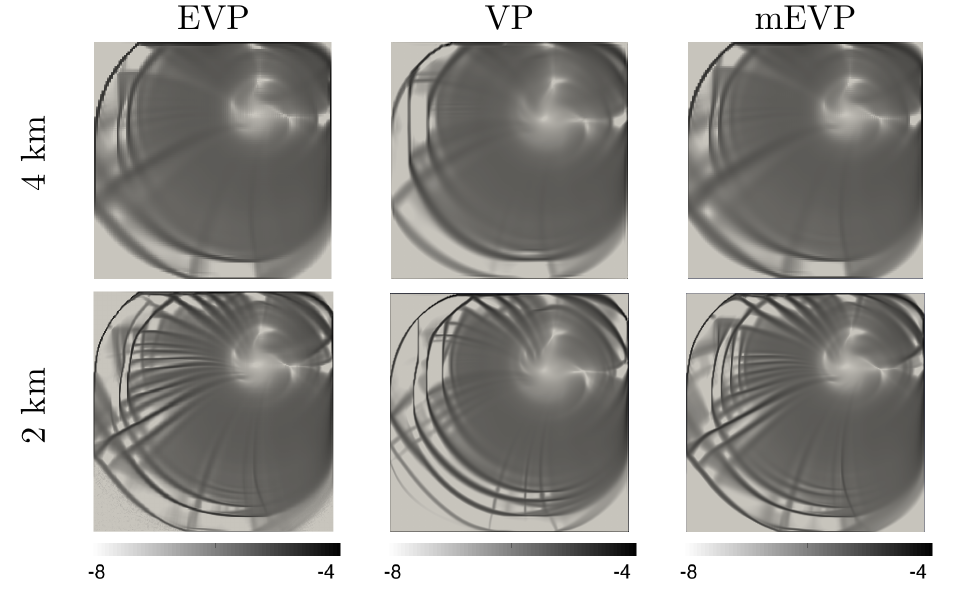}
  \end{center}
  \begin{center}
   \caption{ B-grid that only differs by the (m)(E)VP formulation. 
   The shear deformation is given in $s^{-1}$ and plotted in logarithmic scale.
   \label{fig:Bgrid_EVP}}
  \end{center}
\end{figure}

\begin{figure}
    \centering
    \includegraphics[scale=1.0]{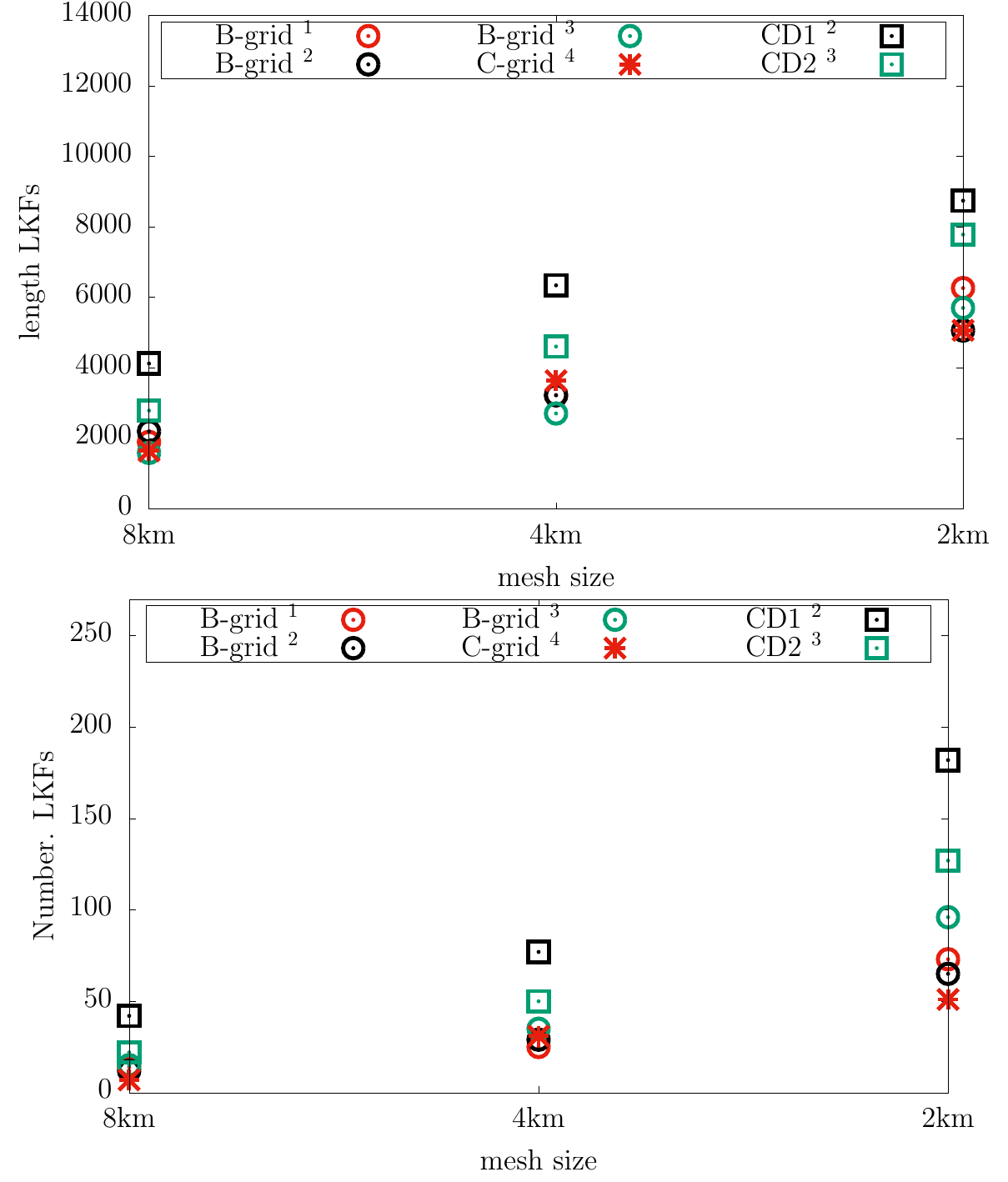}
    \caption{Detected LKFs on quadrilateral grids. The superscripts 1, 2, 3 and 4 in the upper plot refer to the simulations carried out in the framework of CICE, Gascoigne, MPAS and MITgcm respectively.}
   \label{fig:quads_LKFs}
\end{figure}

While Gascoigne solves the VP equations implicitly, MPAS uses the EVP method. We start with pointing out that differences in the solver choice are minor compared to the effect of velocity placement. For this purpose we present simulations based on a B-grid approximation that only differs in the chosen iterative method (see Figure  \ref{fig:Bgrid_EVP}).  All three approach give qualitatively similar results. 
This can be seen by comparing the shear deformation simulated with CD1 (Figure \ref{fig:quads}) to shear deformation presented in Figure \ref{fig:Bgrid_EVP} or by analyzing the amount of detected features in Figure \ref{fig:quads_LKFs}.

The CD1 and CD2 approximation resolve more structure with respect to the cell size ( 8 km, 4 km, 2 km)  than the corresponding B-grid approach in MPAS and Gascoigne (see  Figure \ref{fig:quads} (column 1-3)  and Figure \ref{fig:quads_LKFs}). 
The CD1 approach in Gascoigne and the CD2 approximation in MPAS share the velocity  and tracer placement and therefore have the same number of dof. Apart from that the two setups differs by the used numerical methods to discretize the equations. Figure \ref{fig:quads} shows that the CD2 setup simulates less structure than the CD1 approach. This observations are confirmed by the detection algorithm (see Figure \ref{fig:quads_LKFs}).

\section{Discussion}\label{sec:diss}
Due to uncertainties in the detection algorithm \citep{Mehlmannetal2021} we discuss only results that exhibit a consistent trend across all three measures - the number of LKFs, the length of LKFs, and visual evaluation of the approximation. 

{
Please note that in our study the resolving capacity of the different discretizations is measured in terms of simulated LKFs. It remains to be shown that the simulated deformations reproduces observed scaling characteristics in Pan Arctic sea-ice simulations. (E)VP models reproduce scaling characteristics once the mesh resolution is high enough \citep{Bouchat2022}. As main deformation characteristics are present in all the considered discretizations (see Figure \ref{fig:trishear}, Figure \ref{fig:hexshear} and Figure \ref{fig:quads}) it is likely that the discussed setups will simulate observed scaling characteristics to a certain extent. }

The numerical analysis in Section \ref{sec:num} shows that among the discretizations considered here the CD-grids resolve more LKFs for given grid resolution than the A-grid on triangular meshes, the hexagonal B-grid and the quadrilateral B-grid/C-grid. This can mainly be attributed to the fact that the CD-grids double the dof on quadrilateral grids, triples the dof on triangular meshes compared to A-grids and increase the dof by factor 1.5 compared to hexagonal B-grids. 

The number of dof for a given grid spacing is a key parameter for the simulation of LKFs. Our analysis considers low order approximations. For the benchmark problem a second order conforming finite element discretization has been tested by \cite{Shih2022}. The discretization showed higher resolving capacity than the first order conforming finite element approximation. It would be of interested to compare the performance of high order spatial discretizations to the presented CD-grids.

Overall the nonconforming CD1 approximation produces more deformation structure for given  grid resolution than the CD2 setup on both quadrilateral and hexagonal/triangular meshes, even though the CD1 and CD2 approximation have the same number of dof.

The lower resolution capacity with respect to the mesh size can be explained as follows.
In the quadrilateral case the CD2 approach can be interpreted as a rotated B-grid with a doubled number of cells. For simplicity  we consider a unit square with $\text{N}=\frac{1}{h^2}$ cells. The rotated B-grid can be viewed as a diamond shaped element that is placed in each quadrilateral cell. The side length of the diamond is given by $\frac{h}{\sqrt{2}}$. This means that in case of the rotated grid the number of elements rise to $2\text{N}=\frac{2}{h^2}$. Therefore the rotated B-grid has the same resolution as the B-grid on meshes with a side length of $h \sqrt{2}$ e.g $\frac{N}{2}$ cells.
The CD1 approach simulates the same deformation structure as the B-grid but on meshes twice as coarse \citep{Mehlmannetal2021}. Thus the CD1 grid has the same resolution as the B-grid on meshes with a grid spacing of $2h$ e.g $\frac{N}{4}$. This shows that CD1 approach simulates more structure w.r.t. the mesh size than the CD2 discretization and the B-grid simulates less LKFs than the CD2 framework. 

The resolving capacity on triangular grids is discussed in \cite{Danilov2022} and in \cite{Danilov2023}. Based on a Fourier analysis the authors show that the CD1 approximation provides the highest accuracy, followed by the CD2 approximation, and the A-grid.

 
 
As shown by \cite{Mehlmannetal2021} the nonconforming CD1 approach resolves more LKFs than quadrilateral B-grids/C-grids and triangular A-grids but on meshes with doubled grid spacing (less dof).
This conclusion can not be drawn for the CD2 setup (see Figure \ref{fig:quads_dof_LKFs}). The CD2 approximation does not provide the same resolution capacity per dof neither on quadrilateral nor on triangular/hexagonal meshes.

On hexagonal meshes with the same dof in the velocity the CD2 simulates more LKFs than the B-grid (see Figure{ \ref{fig:quads_dof_LKFs}}). We can not confirm this result on quadrilateral meshes (see Figure upper plot \ref{fig:quads_dof_LKFs}). 
We attribute such a different behavior on quadrilateral and hexagonal meshes to the different ratio of velocity and tracer dof. Note that the comparison on grids with the same number of dof is slightly in favor of the B-grid because doubling the number of cells doubles the dof of the tracers. It has been shown by \cite{MehlmannDanilov2022} that more dof in the tracer point promotes the production of LKFs as the thickness and concentration influences the representation of the pressure $P$ in the rheology (\ref{eq:stress}).
The study of \cite{Danilov2023} compares the triangular CD2 approximation to an A-grid discretization with the same dof and shows that the A-grid produces more LKFs than the CD2 discretization. On meshes with the same velocity dof the tracer points are tripled in the A-grid, which leads to an improved representation of the tracers. 

In case of the benchmark problem, implicit solving VP and explicit subcyling (m)EVP leads to qualitatively similar results (see Figure \ref{fig:Bgrid_EVP}). We conclude that the differences in the simulation of LKFs introduced by varying the solvers are much smaller than differences between CD1 and CD2 grids. The CD1 discretization needs a stabilization. In case of implicit solvers the configuration of the stabilizing parameter is identified \citep{Mehlmann2021, Mehlmannetal2021}, while the optimal choice in the context of (m)EVP methods requires additional tuning. Using the CD2 setup no additional stabilization is needed, which is an advantage of the CD2 approach.



\section{Conclusion}\label{sec:conc}
 We find that the nonconforming CD1 approximation \citep{Mehlmann2021, Mehlmannetal2021} produces more deformation structure than the CD2 approach \citep{Capodaglio2022} on both quadrilateral and triangular meshes, even so 
 both CD-grids have the same number of dof. This shows that besides the placement of the velocity, the chosen spatial discretization plays an important role for the simulation of LKFs.

The nonconforming CD1 approach {provides a promising resolution property.}
Even though CD-grids doubles the number of dof compared to quadrilateral A, B and C-grids and  triples the dof of the triangular A-grid, the CD1 setup simulates qualitatively similar LKFs on meshes with half of the grid spacing (4 times less dof). We found that the CD2 setup does not have this resolution capacity. However on hexagonal meshes with the same number of dof the MPAS CD2 approach simulates more LKFs than the MPAS B-grid setup.

The CD2 discretization resolves more LKFs than standard A,B or C-grids on a fixed mesh. This can mainly be attributed to the higher amount of velocity dof. Nevertheless, the higher number of dof is an appealing property as for the simulation of deformation structure in viscous-plastic sea-ice models a high spatial mesh resolution is needed.

\section{Data and software availability statement}
 The version of MPAS-Seaice used for the results in this paper is provided at  \url{https://doi.org/10.5281/zenodo.7662610} \citep{Giacomodata2023}.
 The sea-ice component of FESOM used for simulations reported here is available at \url{https://doi.org/10.5281/zenodo.7646908} \citep{danilovscode2023}.
The Gascoigne, ICON, MITgcm, CICE, FESOM (A-grid and CD1) data is available at \citep{Mehlmanndata2021}. The MPAS and FESOM CD2 data and the routines to process it can be accessed via \url{ https://data.mendeley.com/datasets/7h9hkjvx48/1 } \citep{Mehlmanndata2023}.

\section*{Acknowledgment}
We thank N. Hutter for the development of his LKF diagnostics.
 C. Mehlmann is funded by the the Deutsche Forschungsgemeinschaft
(DFG, German Research Foundation) - Project number 463061012. Document approved for unlimited release: LA-UR-23-21941. This paper is a contribution to the project S2: Improved parameterisations and numerics in climate models of the Collaborative Research Centre TRR 181 “Energy Transfer in Atmosphere and Ocean” funded by the Deutsche Forschungsgemeinschaft (DFG, German Research Foundation) – project no. 274762653.

\end{document}